\documentclass[review]{siamart1116}
\usepackage{amssymb}
\usepackage{mathtools}
\usepackage{amsmath}
\usepackage{graphicx}
\usepackage{grffile}

\usepackage{lipsum}
\usepackage{amsfonts}
\usepackage{graphicx}
\usepackage{epstopdf}
\usepackage{algorithmic}
\ifpdf
  \DeclareGraphicsExtensions{.eps,.pdf,.png,.jpg}
\else
  \DeclareGraphicsExtensions{.eps}
\fi

\numberwithin{theorem}{section}

\newcommand{\TheTitle}{Chebyshev Interpolation For Function in 1D } 
\newcommand{\TheAuthors}{Tianyu Sun}
\headers{\TheTitle}{\TheAuthors}

\title{\TheTitle}

\author{
Tianyu Sun\\
Advisor: Prof. Mark Holmes
}

\usepackage{amsopn}

\ifpdf
\hypersetup{
  pdftitle={\TheTitle},
  pdfauthor={\TheAuthors}
}
\fi

\begin{document}

\maketitle

\begin{abstract}
  This research is concerned with finding the roots of a function in an interval using Chebyshev Interpolation. Numerical results of Chebyshev Interpolation are presented to show that this is a powerful way to simultaneously calculate all the roots in an interval.
\end{abstract}

\begin{keywords}
  Chebyshev Interpolation, Chebyshev-Frobenius Matrix, Root Finding
\end{keywords}

\section{Introduction}
Classically, people use Newton's iteration or secant method to calculate the roots of a given function\cite{Holmes}. But the problem is that those methods can only calculate one root at a time, and sometimes the method can fail if the initial point is not chosen close to the solution. Here, we investigate a rootfinding method that uses Chebyshev interpolation\cite{Boyd} and explore its capabilities in finding the roots in an interval. \\

Let $f(x)$ be an infinitely differentiable function defined on $[a,b]$. To compute the roots of $f(x)$ over $[a,b]$, the first step is to replace the function by a polynomial approximation $f_N(\hat{x})$ in Chebyshev series form since the key idea is to interpolate the function in Chebyshev coordinate and then transform back to ordinary coordinate. In other words, let 
\[
f(x)\approx f_N(\hat{x})=a_0+a_1T_1(\hat{x})+a_2T_2(\hat{x})+...+a_{N}T_{N}(\hat{x})
\]\\

Each $T_n(\hat{x})$ is the Chebyshev polynomial of the first kind of degree $n$. The next step is to put the coefficient\cite{CI}\cite{CII} of the series $a_1,a_2,...$  into a matrix called Chebyshev-Frobenius matrix\cite{Boyd}. Then we calculate the eigenvalues of the matrix, and they are all possibilities of the roots of the function. Finally, we can use Newton's method to refine and get rid of some spurious roots.\\

\section{Chebyshev Interpolation}
\label{sec:main}

\subsection{Continuous Chebyshev expasion}

A function $f(x)$ can be approximated in terms of Chebyshev polynomials as
\begin{equation}
f(x)\approx f_N(\hat{x})=\sum_{j=0}^{N} a_jT_j(\hat{x})
\end{equation}
where $T_j({\hat{x}})$ is the Chebyshev polynomial of the first kind, N is the number of terms, and $\hat{x}=\frac{2x-(a+b)}{b-a}$.\\

The Chebyshev polynomial of the first kind is defined as 
\begin{equation}
T_j(x)=\cos(j\arccos(x))
\end{equation}
Then we have 
\begin{equation*}
\begin{split}
T_0(x)&=1\\
T_1(x)&=x\\
T_2(x)&=2x^2-1\\
\vdots
\end{split}
\end{equation*}

Chebyshev polynomials are orthogonal and enjoy the following orthogonality relation:
\begin{equation*}
\int_{-1}^{1}T_j(x)T_k(x)(1-x^2)^{-1/2}dx=N_{jk}\delta_{jk}
\end{equation*}
with $N_{00}=\pi$ and $N_{rr}=\frac{1}{2}\pi$ if $r\not= 0$. $\delta_{jk}$ is the Kronecker Delta function and defined as: $\delta_{jk}=1$ if $j=k$ and $\delta_{jk}=0$ if $j\not =k$

Multiplying the coefficient $a_j$ on both sides, we get
\begin{equation*}
\begin{split}
\int_{-1}^{1}a_jT_j(x)T_k(x)(1-x^2)^{-1/2}dx&=a_jN_{jk}\delta_{jk}\\
\int_{-1}^{1}f(x)T_k(x)(1-x^2)^{-1/2}dx&=a_jN_{jk}\delta_{jk}\qquad \left(f(x)=a_jT_j(x)\right)\\
\int_{-1}^{1}\frac{f(x)T_j(x)}{\sqrt{1-x^2}}dx&=a_jN_{jj}\qquad \left(\delta_{jk}=0 \mbox{ if } j\not = k\right)
\end{split}
\end{equation*}
\begin{equation}
\begin{split}
a_0&=\pi\int_{-1}^{1}\frac{f(x)}{\sqrt{1-x^2}}dx\\
a_j&=\frac{\pi}{2}\int_{-1}^{1}\frac{f(x)T_j(x)}{\sqrt{1-x^2}}dx \qquad\qquad(\mbox{for  } j\not= 0)
\end{split}
\end{equation}

\subsection{Discrete Chebyshev expasion}

When the integral can not be evaluated exactly, we can introduce a discrete grid and use a numerical formula. Chebyshev polynomials are orthogonal and enjoy the following discrete orthogonality relation:

For $i,j\leq N$,
\begin{equation}
\sum_{k=0}^N T_i(\hat{x_k})T_j(\hat{x_k})=
\begin{cases}
0 & i\not= j\\
N/2 & i=j\not=0\\
N & i=j=0
\end{cases}
\end{equation}

\noindent where $\hat{x_k}$ is the zero of Chebyshev polynomial at N-th degree and is defined as 
\begin{equation*}
\hat{x_k}=\cos\frac{\pi(2k-1)}{2N}
\end{equation*}

Since $f_N(\hat{x_k})$ interpolates $f$ at the $N+1$ Chebyshev nodes, we have that at these nodes $f(x_k)=f_N(\hat{x_k})$. Then we have
\begin{equation}
\sum_{k=0}^N f(x_k)T_i(\hat{x_k})=\sum_{j=0}^{N} a_j\left[\sum_{k=0}^N T_i(\hat{x_k})T_j(\hat{x_k})\right]
\end{equation}

Using (4), we can find that 
\begin{equation}
\begin{split}
a_0&=\frac{1}{N}\sum_{k=1}^N f(x_k)\\
a_j&=\frac{2}{N}\sum_{k=1}^N f(x_k)T_j(\hat{x_k})\qquad\qquad(\mbox{for  } j\not= 0)
\end{split}
\end{equation}
where $x_k$ is the zero of Chebyshev polynomia in the whole interval and can be calculated using $x_k=\frac{1}{2}(a+b+(b-a)\hat{x_k})$. $f(x_k)$ is the value of original function evaluated at these zeros. \\

$T_j(x_k)$ is the value of Chebyshev polynomial $T_j(x)$ at the zeros of $T_N(x)$. Using equation (2), $T_j(x_k)$ can be found as
\begin{equation*}
T_j(x_k)=\cos\frac{j(2k-1)\pi}{2N}
\end{equation*}\\

\subsection{Rate of Convergence}

If we repeatedly integrate equation (3) by parts, we get 

\begin{equation*}
a_n=\frac{1}{n^m}\frac{2}{\pi}\int_{-1}^{1} \frac{f^{(m)}(x)T_j(x)}{\sqrt{1-x^2}} dx
\end{equation*}

\noindent Thus, if $f$ is m-times differentiable in $[-1,1]$, the above integral exists, and we can conclude that $a_n=O(n^{-m}), n=1,2,...$.\\

Since $T_j$ is bounded above by $1$ on $[-1,1]$, it follows that the truncation error is bounded by the sum of the absolute value of negelected coefficients:

\begin{equation*}
\mid f(x)-f_N(X)\mid \leq\sum_{n=N+1}^{\infty} \mid a_n\mid
\end{equation*}

\noindent Therefore, we can conclude that when a function $f$ has $m + 1$ continuous derivatives on $[−1, 1]$,where $m$ is a finite number, then $\mid f(x)-f_N(x)\mid = O(n^{-m})$ as $n\to \infty$ for all $x\in [−1, 1]$. If $f$ is infinitely differentiable, then the convergence is faster than $O(n^{-m})$ no matter how large we take $m$. 

\section{Chebyshev-Frobenius matrix}
\label{sec:alg}

From equation (1), we get a polynomial expansion $f_N(\hat{x})=\sum_{j=0}^N a_jT_j(\hat{x})$ for $f(x)$. If we multiply each basis function $T_j(\hat{x})$ by $\hat{x}$, we get a new polynomial of degree $j+1$. We can reexpand the polynomial as
\begin{equation*}
\hat{x}T_j(\hat{x})=\sum_{k=0}^{j+1} H_{j,k}T_j(\hat{x})
\end{equation*}
for some coefficients $H_{j,k}$.

Define a vector $\vec{T}$ whose $N$ elements are the basis functions. Then the first $N$ polynomials $xT_j(\hat{x})$ can be organized in to a matrix equation,
\begin{equation*}
{\bf H}\vec{T}=\hat{x}\vec{T}-H_{N,N+1}T_{N}(\hat{x})\vec{e_{N}}
\end{equation*}
Note that the term $H_{N+1}T_{N}(\hat{x})\vec{e_N}$ arises because the product of the last term $xT_{N-1}$ is a polynomial of degree $N$. However, $T_{N}\not\in \vec{T}$. Therefore, we have to substract the term to make both sides of the equation balance. \\

We can remove the extra $T_{N}$ term by adding $qf_N(\hat{x})$ and substracting it from the last row of the matrix, which becomes
\begin{equation*}
\begin{split}
&\sum_{k=0}^{N-1}H_{N,k+1}T_k(\hat{x})\\
=&\hat{x}T_{N-1}(\hat{x})-H_{N,N+1}T_{N}(\hat{x})+q\left\{a_NT_N(\hat{x})+\sum_{j=0}^{N-1}a_jT_j(\hat{x})\right\}-qf_N(\hat{x})
\end{split}
\end{equation*}

If $q=H_{N,N+1}/a_N$, the $T_{N}$ term can be cancelled. Then the equation becomes
\begin{equation*}
\begin{split}
\sum_{k=0}^{N-1}H_{N,k+1}T_k(\hat{x})=\hat{x}T_{N-1}(\hat{x})+\frac{H_{N,N+1}}{a_N}\sum_{j=0}^{N-1}a_jT_j(\hat{x})-\frac{H_{N,N+1}}{a_N}f_N(\hat{x})\\
\sum_{k=0}^{N-1}\left\{H_{N,k+1}T_k(\hat{x})-H_{N,N+1}\frac{a_k}{a_N}T_k(\hat{x})\right\}=\hat{x}T_{N-1}-\frac{H_{N,N+1}}{a_N}f_N(\hat{x})
\end{split}
\end{equation*}

If $\hat{x}$ is the root of $f_N(\hat{x})$, so $f_N(\hat{x})=0$. Then the matrix problem becomes an eigenvalue problem, ${\bf H}\vec{T}=\hat{x}\vec{T}$ where the elements of ${\bf M}$ are
\begin{equation*}
\begin{split}
M_{j,k}=H_{j,k}&\quad j=1,2,...N-1\\
M_{N,k}=H_{N,k}-H_{N,N+1}\frac{a_{k-1}}{a_N}&\quad k=1,2,...N
\end{split}
\end{equation*}

Based on this, we can find the root using a special matrix called a Chebyshev-Frobenius matrix. This is defined as
\begin{equation*}
A_{jk}=\begin{cases}
\delta_{2,k} &\quad j=1,k=1,2,...,N\\
\frac{1}{2}\{\delta_{j,k+1}+\delta_{j,k-1}\}&\quad j=2,...,(N-1),k=1,2,...,N\\
(-1)\frac{a_{k-1}}{2a_N}+\frac{1}{2}\delta_{k,N-1}\qquad&\quad j=N,k=1,2,...,N
\end{cases}
\end{equation*}
where $\delta_{jk}$ is the usual Kronecker Delta function.

For example, when the polynomial is 5-th degree, we can write the matrix as
\begin{equation*}
\begin{bmatrix}
0 & 1 & 0 & 0 & 0\\
1/2 & 0 & 1/2 & 0 & 0\\
0 & 1/2 & 0 & 1/2 & 0\\
0 & 1 & 1/2 & 0 & 1/2\\
(-1)\frac{a_0}{2a_5} & (-1)\frac{a_1}{2a_5} & (-1)\frac{a_2}{2a_5} & (-1)\frac{a_3}{2a_5}+(1/2) & (-1)\frac{a_4}{2a_5}\\
\end{bmatrix}
\end{equation*}

The eigenvalues of the matrix are roots of Chebyshev nodes in standard interval. We only accept roots that are within $10^{-8}$ of real axis and within $10^{-6}$ of 1.\\

Finally, we can convert the roots from standard interval to the whole interval using $x=\frac{1}{2}[a+b+(b-a)\hat{x}]$ and get all the roots in $x\in[a,b]$.\\

\section{Newton-Polishing}

After calculating all possible roots, we can examine the validity of each root locally by using one or two Newton iterations of $f(x)$. Newton's iteration is defined as 
\begin{equation*}
x^{(n+1)}=x^{(n)}-\frac{f(x^{(n)})}{df/dx(x^{(n)})}
\end{equation*}

Therefore, once we find a possible root, we can refine the root by finding a better root close to it. For efficiency, it is important to stop when further iterations produce no reduction in the correction. \\

\section{Numerical Experiment}
\label{sec:experiments}

Chebyshev Rootfinder works very well for various continuous functions. First we try to show how well Chebyshev approximation can interpolate a certain function by comparing the graph of actual function and polynomial obtained by certain N-th degree Chebyshev Interpolation. Then we try to show how well Chebyshev Interpolation converges by comparing the roots obtained from different N-th degree Chebyshev Interpolation.

\subsection{Case 1: Trigonometrix function}
For the first case, we use trigonometrix function
\begin{equation*}
f(x)=\cos(x)
\end{equation*}

where $x\in[-10,10]$. We compare 12-th degree Chebyshev Interpolation with the real function since the difference becomes so small if we take more terms to approximate the function. We also find all possible roots obtained using Chebyshev-Frobenius matrix and all accepted roots. We only accept roots that are within $10^{-8}$ of real axis and within $10^{-6}$ of 1. Finally, we convert those roots to real axis and see how those roots converge as we take more term to approximate the function. The results are illustrated in Figure 1,2,3,4, and 5.\\
\includegraphics[width=\textwidth]{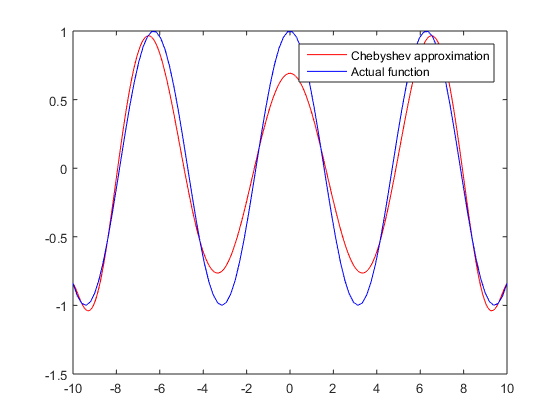}
\begin{center}
Figure 1: Graph of actual function and Chebyshev approximation for $N=12$. \\
\end{center}
\includegraphics[width=\textwidth]{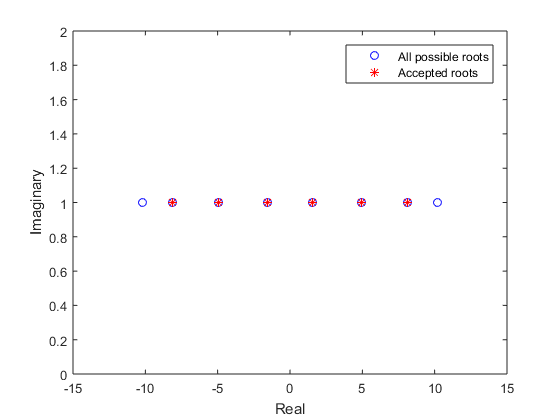}
\begin{center}
Figure 2: All possible roots for $N=13$. \\
\end{center}
\includegraphics[width=\textwidth]{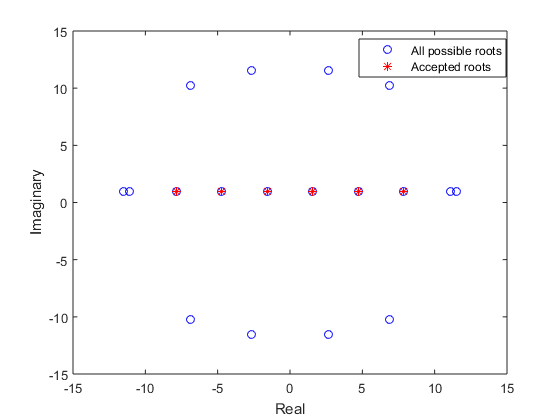}
\begin{center}
Figure 3: All possible roots for $N=20$. \\
\end{center}
\includegraphics[width=\textwidth]{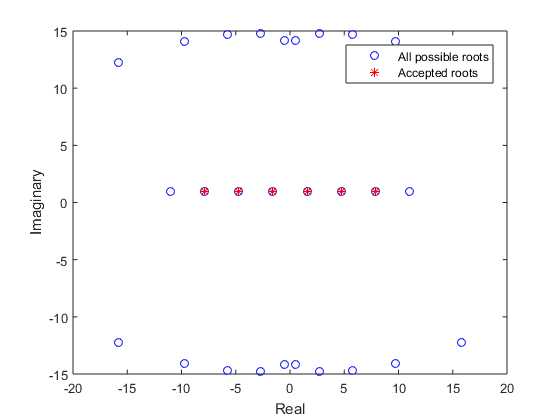}
\begin{center}
Figure 4: All possible roots for $N=30$. \\
\end{center}
\includegraphics[width=\textwidth]{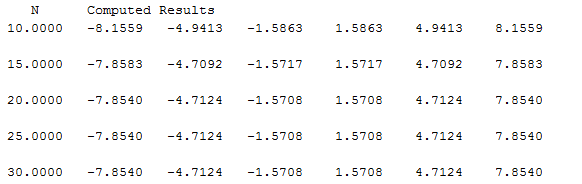}
\begin{center}
Figure 5: All acceped roots in $x\in[-10,10]$ for different value of $N$.\\
\end{center}

As can be seen, when we take more Chebyshev terms to approximate the function, the number of possible roots increases. Since we only accept roots in certain interval, the number of accepted roots stays the same. Also, notice that the accuracy of the actual results increases as we take more and more Chebyshev terms, which suggests the convergence of Chebyshev approximation.

\subsection{Case 2: Exponential function}
For the second case, we use exponential function
\begin{equation*}
f(x)=e^x
\end{equation*}

where $x\in[-10,10]$. We compare 8-th degree Chebyshev Interpolation with the real function for the same reason as the previous one. We also find all possible roots obtained using Chebyshev-Frobenius matrix and all accepted roots.  Finally, we convert those roots to real axis and see how those roots converge as we take more term to approximate the function. The results are illustrated in Figure 6,7,8,9, and 10.\\
\includegraphics[width=\textwidth]{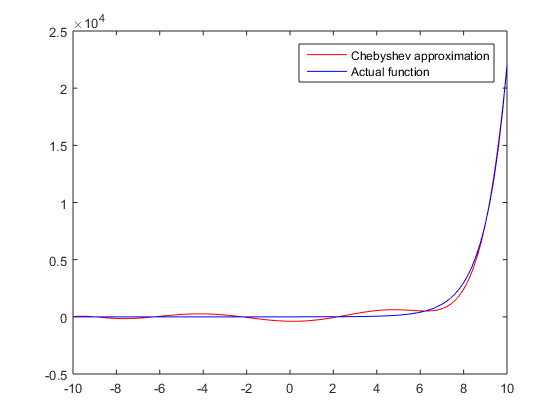}
\begin{center}
Figure 6: Graph of actual function and Chebyshev approximation for $N=8$. \\
\end{center}
\includegraphics[width=\textwidth]{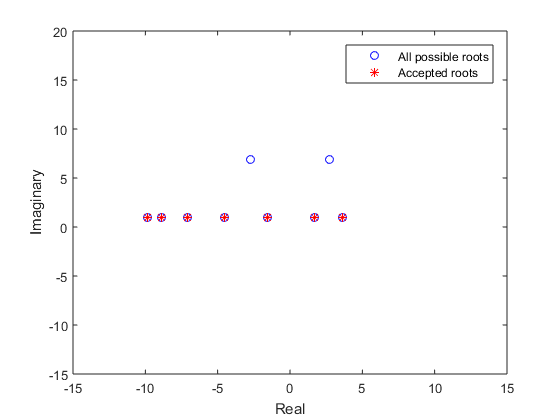}
\begin{center}
Figure 7: All possible roots for $N=13$. \\
\end{center}
\includegraphics[width=\textwidth]{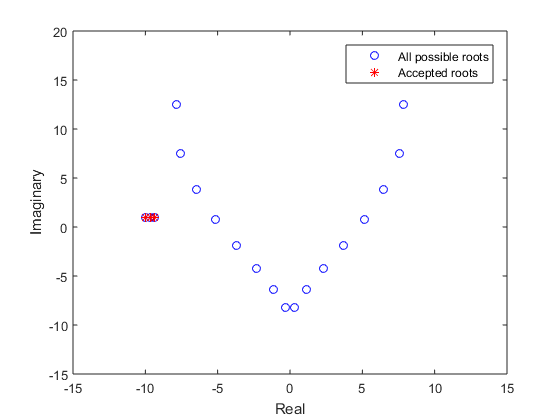}
\begin{center}
Figure 8: All possible roots for $N=20$. \\
\end{center}
\includegraphics[width=\textwidth]{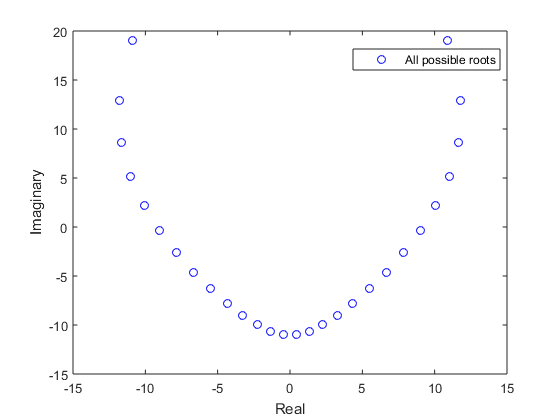}
\begin{center}
Figure 9: All possible roots for $N=30$. \\
\end{center}
\includegraphics[width=\textwidth]{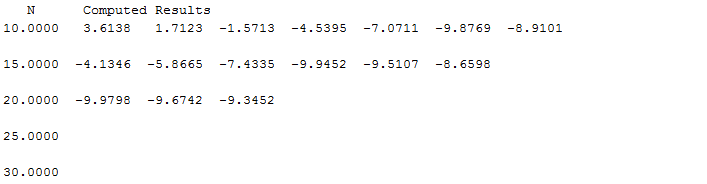}
\begin{center}
Figure 10: All acceped roots in $x\in[-10,10]$ for different value of $N$.\\
\end{center}

As can be seen, when we only take only few terms to approximate the function, it gave some spurious roots. However, when we take more and more terms to approximate the function, the number of spurious and eventually disappears, which also suggests the convergence of Chebyshev approximation.

\subsection{Case 3: Slightly complex function}
Chebyshev approximation doesn't work very well some slightly complex function in that the number of calculated roots may increase as we take more and more terms. For the last case, we use 
\begin{equation*}
f(x)=e^{(-0.5x^2)}(12-48x^2+16x^4)
\end{equation*}

where $x\in[-10,10]$. We compare 30-th degree Chebyshev Interpolation with the real function to understand why the number of calculated roots increase as the term of Chebyshev Interpolation increases. We also find all possible roots obtained using Chebyshev-Frobenius matrix and all accepted roots.  Finally, we convert those roots to real axis to show that the number of roots increases as the number of term in Chebyshev Interpolation increases. The results are illustrated in Figure 11,12,13,14, and 15.\\
\includegraphics[width=\textwidth]{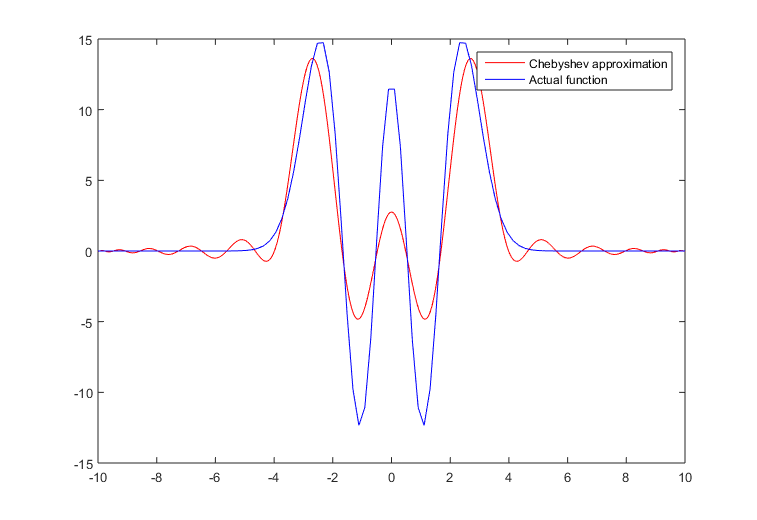}
\begin{center}
Figure 11: Graph of actual function and Chebyshev approximation for $N=30$.\\
\end{center}
\includegraphics[width=\textwidth]{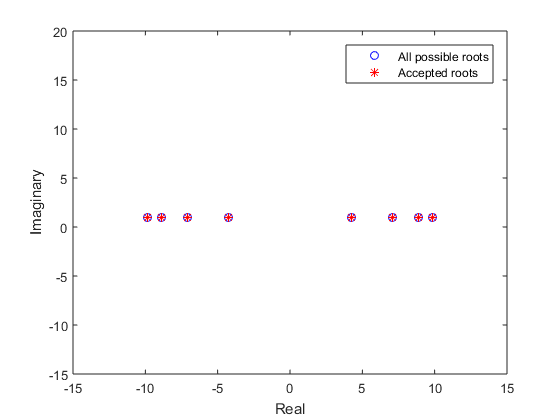}
\begin{center}
Figure 12: All possible roots for $N=10$. \\
\end{center}
\includegraphics[width=\textwidth]{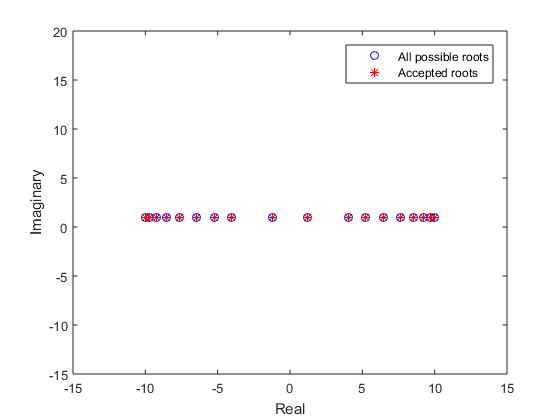}
\begin{center}
Figure 13: All possible roots for $N=20$. \\
\end{center}
\includegraphics[width=\textwidth]{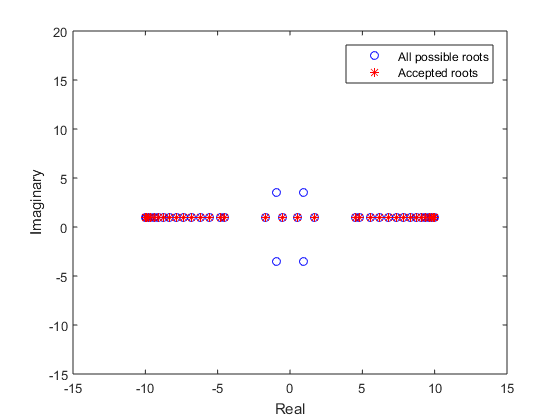}
\begin{center}
Figure 14: All possible roots for $N=40$. \\
\end{center}
\includegraphics[width=\textwidth]{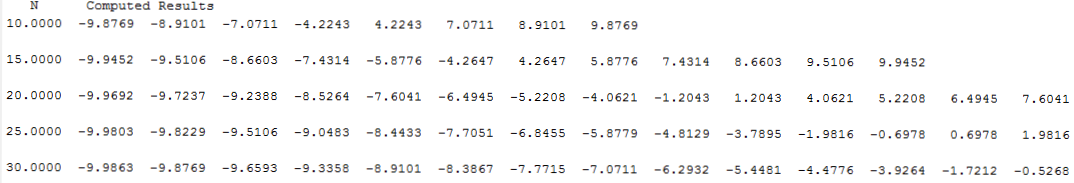}
\begin{center}
Figure 15: All acceped roots in $x\in[-10,10]$ for different value of $N$.\\
\end{center}

From the figure 11, we can see that there are oscillations near the end of boundaries. This phenomenon is know as Gibb's phenomenon. Since Chebyshev approximation is a great tool to find all possible roots globally, we can use other interpolation method like Newton's method to get rid of those spurious roots and find the actual roots locally.

\section{Conclusions}
\label{sec:conclusions}

As one can see, Chebyshev approximation has a really fast convergence rate. Unlike Newton's method, it is a great tool to compute all possible roots at once globally. However, sometimes Chebyshev approximation can yield spurious roots due to Gibb's phenomenon, and it is a good idea to use Newton's method in the end to refine individual roots locally.

\bibliographystyle{siamplain}
\bibliography{references}
\end{document}